\documentclass{elsart3-1}
\usepackage[english,francais]{babel}   
\usepackage{amsmath}
\usepackage{amsfonts}
\usepackage{amssymb}
\usepackage{bbm}

\newtheorem{theorem}{Theorem}[section]
\newtheorem{lemma}[theorem]{Lemma}
\newtheorem{proposition}[theorem]{Proposition}
\newtheorem{corollary}[theorem]{Corollary}
\newtheorem{definition}[theorem]{Definition\rm}

\newtheorem{remark}{\it Remark\/}

\newtheorem{theoreme}{Th\'eor\`eme}[section]

\newtheorem{f-proposition}[theoreme]{Proposition}

\newtheorem{f-definition}[theoreme]{D\'efinition\rm}

\def\og{\leavevmode\raise.3ex\hbox{$\scriptscriptstyle\langle\!\langle$~}}
\def\fg{\leavevmode\raise.3ex\hbox{~$\!\scriptscriptstyle\,\rangle\!\rangle$}}

\newcommand{\proof}{\paragraph*{Proof}}
\newcommand{\fin}{$\Box$\\}

\newcommand{\sifbm}{\mathbf{B}}

\begin{document}

\centerline{}
\begin{frontmatter}

\selectlanguage{english}
\title{A Characterization of the Set-indexed Fractional Brownian Motion by Increasing Paths}

\author[author1]{Erick Herbin}\ead{erick.herbin@dassault-aviation.fr} ,
\author[author2]{Ely Merzbach}\ead{merzbach@macs.biu.ac.il}
\address[author1]{Dassault Aviation, 78 quai Marcel Dassault, 92552 Saint-Cloud Cedex,
France} 
\address[author2]{Dept. of Mathematics,
Bar Ilan University, 52900 Ramat-Gan, Israel}

\begin{abstract}
\selectlanguage{english}
We prove that a set-indexed process is a set-indexed fractional Brownian motion
if and only if its projections on all the increasing paths are one-parameter
time changed fractional Brownian motions. As an application, we present an
integral representation for such processes.

\vskip 0.5\baselineskip

\selectlanguage{francais}
\noindent{\bf R\'esum\'e} \vskip 0.5\baselineskip
\noindent
{\bf Une caract\'erisation par chemins croissants du mouvement brownien fractionnaire index\'e par des ensembles.}
On montre qu'un processus stochastique est un mouvement brownien fractionnaire index\'e par des ensembles si et seulement si ses projections sur tous les chemins croissants sont des mouvements browniens fractionnaires \`a param\`etres r\'eels chang\'es de temps. On applique ce r\'esultat \`a la d\'efinition d'une repr\'esentation int\'egrale pour de tels processus.

\end{abstract}

\end{frontmatter}

\selectlanguage{francais}
\section*{Version fran\c{c}aise abr\'eg\'ee}

Dans \cite{ehem}, le {\em mouvement brownien fractionnaire index\'e par des ensembles (sifBm)} est d\'efini et ses propri\'et\'es de stationnarit\'e et d'autosimilarit\'e sont \'etudi\'ees. D'autre part, on prouve que la projection d'un sifBm sur un chemin croissant est un mouvement brownien fractionnaire index\'e par $\mathbf{R}_+$ chang\'e de temps. L'objet de cette note est la r\'eciproque de ce r\'esultat.

En consid\'erant une collection d'indices $\mathcal{A}$ de sous-ensembles compacts d'un espace m\'etrique localement compact $\mathcal{T}$ muni d'une mesure de Radon $m$ (voir \cite{ehem}), le  mouvement brownien fractionnaire index\'e par $\mathcal{A}$ est d\'efini comme le processus gaussien centr\'e $\sifbm^H=\left\{\sifbm^H_U;\;U\in\mathcal{A}\right\}$ tel que
\begin{equation}\label{defsifbm}
\forall \ U,V\in\mathcal{A};\  E\left[\sifbm^H_U \sifbm^H_V\right]=\frac{1}{2}
\left[m(U)^{2H}+m(V)^{2H}-m(U\bigtriangleup V)^{2H}\right],
\end{equation}
o\`u $0<H\leq\frac{1}{2}$.

\begin{f-definition}
On appelle {\em flot \'el\'ementaire} toute fonction $f:[a,b]\subset\mathbf{R}_+\rightarrow\mathcal{A}$ v\'erifiant
\begin{align*}
\forall s,t\in [a,b];\quad & s<t \Rightarrow f(s)\subseteq f(t)\\
\forall s\in [a,b);\quad & f(s)=\bigcap_{v>s}f(v)\\
\forall s\in (a,b);\quad & f(s)=\overline{\bigcup_{u<s}f(u)}.
\end{align*}
\end{f-definition}

\begin{f-definition}
Un processus index\'e par des ensembles $X=\left\{X_U;\;U\in\mathcal{A}\right\}$ est dit {\em continu monotone ext\'erieurement dans $L^2$} si $X_U$ est carr\'e int\'egrable pour tout $U\in\mathcal{A}$ et
pour toute suite d\'ecroissante $\left(U_n\right)_{n\in\mathbf{N}}$ d'ensembles dans $\mathcal{A}$, $$E\left[|X_{U_n}-X_{\bigcap_m U_m}|^2\right]\rightarrow 0$$ 
quand $n\rightarrow\infty$.
\end{f-definition}

\begin{theoreme}
Soit $X=\left\{ X_U ;\;U\in\mathcal{A} \right\}$ un processus continu monotone ext\'erieurement dans $L^2$.

Si la projection $X^f$ de $X$ sur tout flot \'el\'ementaire $f$, est, \`a un changement de temps pr\`es, un mouvement brownien fractionnaire index\'e par $\mathbf{R}_+$ de param\`etre $H\in (0,1/2)$,
alors $X$ est un mouvement brownien fractionnaire index\'e par $\mathcal{A}$.
\end{theoreme}

Cette caract\'erisation fournit une bonne justification de la d\'efinition du sifBm et ouvre la porte \`a une grande vari\'et\'e d'applications. La repr\'esentation int\'egrale (\ref{intrep}) constitue l'une d'entre elles.

\selectlanguage{english}
\section{Introduction}
In \cite{ehem}, the set-indexed fractional Brownian motion (sifBm)
is defined and its properties of stationarity and self-similarity are
discussed. In particular, it is proved that the projection of a sifBm on an
increasing path is a one-parameter time changed fractional motion. In this
note, we prove the converse.

This characterization gives a good justification of the definition of the sifBm
and opens the door to a variety of applications. Here we present one of them:
an integral representation for the sifBm.

We follow \cite{ehem} for the framework and notation. Our processes are indexed
by an indexing collection $\mathcal A$ of compact subsets of a locally compact
metric space $\mathcal{T}$ equipped with a Radon measure $m$.

The {\em set-indexed fractional Brownian motion (sifBm)} was defined as the
centered Gaussian process
$\sifbm^H=\left\{\sifbm^H_U;\;U\in\mathcal{A}\right\}$ such that
\begin{equation}\label{defsifbm}
\forall \ U,V\in\mathcal{A};\  E\left[\sifbm^H_U \sifbm^H_V\right]=\frac{1}{2}
\left[m(U)^{2H}+m(V)^{2H}-m(U\bigtriangleup V)^{2H}\right],
\end{equation}
where $0<H\leq\frac{1}{2}$.

If $\mathcal{A}$ is provided with a structure of group on $\mathcal{T}$, properties of stationarity and self-similarity are studied in \cite{ehem}.
In the special case of $\mathcal{A}=\left\{ [0,t];\;t\in\mathbf{R}^N_+ \right\} \cup\left\{\emptyset\right\}$, we get a multiparameter process called Multiparameter fractional Brownian motion (MpfBm), whose properties are studied in \cite{mpfbm}.

\section{Projection of the sifBm on flows}
The notion of flow is the key to reduce the proof of many theorems. It was
extensively studied in \cite{cime} and \cite{Ivanoff}.

Let $\mathcal A(u)$ denotes the class of finite unions from sets belonging to $\mathcal A$.

\begin{definition}
An {\em elementary flow} is defined to be a continuous increasing function $f:[a,b]\subset\mathbf{R}_+\rightarrow\mathcal{A}$, i. e. such that
\begin{align*}
\forall s,t\in [a,b];\quad & s<t \Rightarrow f(s)\subseteq f(t)\\
\forall s\in [a,b);\quad & f(s)=\bigcap_{v>s}f(v)\\
\forall s\in (a,b);\quad & f(s)=\overline{\bigcup_{u<s}f(u)}.
\end{align*}

A {\em simple flow} is a continuous function $f:[a,b]\rightarrow\mathcal{A}(u)$ such that there exists a finite sequence $(t_0,t_1,\dots,t_n)$ with $a=t_0<t_1<\dots<t_n=b$ and elementary flows $f_i:[t_{i-1},t_i]\rightarrow\mathcal{A}$ ($i=1,\dots,n$) such that
\begin{equation*}
\forall s\in [t_{i-1},t_i];\quad
f(s)=f_i(s)\cup \bigcup_{j=1}^{i-1}f_j(t_j).
\end{equation*}
The set of all simple (resp. elementary) flows is denoted $S(\mathcal{A})$ (resp. $S^e(\mathcal{A})$).
\end{definition}

\begin{proposition}[\cite{ehem}]]\label{prop1}
Let $\sifbm^H$ be a sifBm and $f$ be an elementary flow. Then the process
$(\sifbm^H)^f=\{\sifbm_{f(t)}^H,\ t\in [a,b]\}$ is a time changed fractional
Brownian motion.
\end{proposition}

The aim of this note is to prove the converse to Proposition \ref{prop1}. For
this purpose, we will use the following lemma proved in \cite{cime}.

\begin{lemma}\label{lemprojflow}
The finite dimensional distributions of an additive $\mathcal
A$-indexed process $X$ determine and are determined by the finite dimensional
distributions of the class $\{X^f,\ f\in S(\mathcal A)\}$.
\end{lemma}

\section{Characterisation of the sifBm}
The converse to Proposition \ref{prop1} in the case of $L^2$-monotone outer-continuous set-indexed processes, gives a characterization of the sifBm by its projection on elementary flows.

Recall the following definition (see \cite{Ivanoff})
\begin{definition}
A set-indexed process $X=\left\{X_U;\;U\in\mathcal{A}\right\}$ is said {\em $L^2$-monotone outer-continuous} if $X_U$ is square integrable for all $U\in\mathcal{A}$ and
for any decreasing sequence $\left(U_n\right)_{n\in\mathbf{N}}$ of sets in $\mathcal{A}$, $$E\left[|X_{U_n}-X_{\bigcap_m U_m}|^2\right]\rightarrow 0$$ 
as $n\rightarrow\infty$.
\end{definition}

\begin{theorem}
Let $X=\left\{ X_U ;\;U\in\mathcal{A} \right\}$ be a $L^2$-monotone outer-continuous set-indexed process.

If the projection $X^f$ of $X$ on any elementary flow $f$, is a time-changed one-parameter fractional Brownian motion of parameter $H\in (0,1/2)$,
then $X$ is a set-indexed fractional Brownian motion.
\end{theorem}

\proof
Let $f:[a,b]\rightarrow\mathcal{A}$ be an elementary flow. As the projected process $X^f$ is a time-changed fBm of parameter $H$, we have
\begin{equation}\label{flowfbm}
\forall s,t\in [a,b];\quad
E\left[X^f_t - X^f_s\right]^2=|\theta_f(t) - \theta_f(s)|^{2H}
\end{equation}
where $X^f_t=X_{f(t)}$ and $\theta_f$ is an increasing function.\\
The idea of the proof is the construction of a measure $m$ such that for any $f\in S^e(\mathcal{A})$,
\begin{equation*}
\forall t\in [a,b];\quad \theta_f(t)=m\left[f(t)\right].
\end{equation*}
\vspace{10pt}

For all $U\in\mathcal{A}$, let us define
\begin{equation*}
F^e_U = \left\{ f \in S^e(\mathcal{A}) :\; \exists u_f\in [a,b]; U= f(u_f) \right\}.
\end{equation*} 
As for all $f$ and $g$ in $F^e_U$, 
$\theta_f(u_f)^{2H} = \theta_g(u_g)^{2H} = E\left[X_U\right]^2$,
one can define
\begin{equation}\label{defpsi}
\psi(U)=\theta_f(u_f)=\left(E\left[X_U\right]^2\right)^{\frac{1}{2H}}.
\end{equation}
For all $U$ and $V$ in $\mathcal{A}$ with $U\subset V$, there exists an elementary flow $f$ such that
\begin{equation*}
\exists u_f,v_f\in [a,b];\; u_f \leq v_f;\quad
U=f(u_f) \subset f(v_f)=V
\end{equation*}
Then, as the time-change $\theta_f$ is increasing, $\psi$ is non-decreasing in $\mathcal{A}$.

The definition of $\psi$ on $\mathcal{A}$ can be extended on the collection $\mathcal{C}$ of sets on the form $C=U\setminus\bigcup_{1\leq i \leq n}U_i$ where $U, U_1,\dots,U_n\in\mathcal{A}$, by the inclusion-exclusion formula
\begin{align}\label{psiC}
\psi(C)=\psi(U)-\sum_{i=1}^n \psi\left(U\cap U_i\right)
&+\sum_{i<j}\psi\left(U\cap (U_i\cap U_j)\right) \nonumber\\
&-\dots+(-1)^n \psi\left(U\cap \left(\bigcap_{1\leq i\leq n}U_i\right)\right)
\end{align} 
The definition (\ref{psiC}) of $\psi$ can be easily extended to the set $\mathcal{C}(u)$ of finite unions of elements of $\mathcal{C}$ in the same way.
Then, for all $C_1, C_2\in\mathcal{C}(u)$ such that $C=C_1\cup C_2 \in\mathcal{C}$, 
\begin{equation}\label{psiunionC}
\psi(C_1\cup C_2) = \psi(C_1) + \psi(C_2) - \psi(C_1\cap C_2).
\end{equation}

From the pre-measure $\psi$ defined on $\mathcal{C}$, the function
\begin{equation}\label{defm}
m : E\subset\mathcal{T} \mapsto \inf_{C_i\in\mathcal{C}\atop E\subset\cup C_i}
\sum_{i=1}^{\infty}\psi(C_i)
\end{equation}
defines an outer measure on $\mathcal{T}$ (see \cite{rogers} pp. 9--26). 
Let us show that $m$ defines a Borel measure on the topological space $\mathcal{T}$.

Let $\mathcal{M}_m$ be the $\sigma$-field of $m$-measurable subsets of $\mathcal{T}$. It is known that $m$ is a measure on $\mathcal{M}_m$ (see \cite{rogers}, thm. 3). By definition, any $U\in\mathcal{A}$ is $m$-measurable if 
\begin{equation*}
\forall A\subset U, \forall B\subset\mathcal{T}\setminus U;\quad
m(A\cup B)=m(A)+m(B).
\end{equation*}
As the inequality $m(A\cup B)\leq m(A)+m(B)$ follows from definition of any outer-measure, it remains to show the converse inequality.

Consider any sequence $\left(C_i\right)_{i\in\mathbf{N}}$ in $\mathcal{C}$ such that $A\cup B\subset\bigcup_{i}C_i$. 
The sequence $\left(C_i\right)_{i\in\mathbf{N}}$ can be decomposed by the elements $C_i,\ i\in I$ such that $C_i\cap U=\emptyset$ and the $C_i,\ i\in J$ such that $C_i\subset U$ (if $C_i\cap U\ne\emptyset$ and $C_i\not\subset U$, cut $C_i=C'_i\cup C''_i$ where $C'_i\subset U$ and $C''_i\cap U=\emptyset$).

From
\begin{equation*}
A\cup B \subset\left[\bigcup_{i\in I}C_i\right] \cup 
\left[\bigcup_{i\in J}C_i\right],
\end{equation*}
we get
\begin{align*}
\sum_{i=1}^{\infty}\psi(C_i) = 
\underbrace{\sum_{i\in I}\psi(C_i)}_{\geq m(B)}+
\underbrace{\sum_{i\in J}\psi(C_i)}_{\geq m(A)}
\end{align*}
which leads to $m(A\cup B)\geq m(A)+m(B)$.

We have proved that $\mathcal{A}\subset\mathcal{M}_m$. By definition of $\mathcal{A}$, the smallest $\sigma$-field containing $\mathcal{A}$ is the Borel $\sigma$-field $\mathcal{B}$. Therefore, $\mathcal{B}\subset\mathcal{M}_m$ and $m$ is a measure on $\mathcal{B}$.
\vspace{10pt}

The second part of the proof is to show that the measure $m$ is an extension of $\psi$, i.e. 
\begin{equation}\label{ext}
\forall U\in\mathcal{A};\quad m(U)=\psi(U).
\end{equation}

\begin{itemize}

\item For any $U\in\mathcal{A}$, by definition of $m(U)$,
\begin{equation}\label{ext1}
m(U)=\inf_{C_i\in\mathcal{C}\atop U\subset\cup C_i}
\sum_{i=1}^{\infty}\psi(C_i)
\leq \psi(U).
\end{equation}

\item To prove the converse inequality, consider $U\in\mathcal{A}$ and a sequence $\left(C_i\right)_{i\in\mathbf{N}}$ in $\mathcal{C}$ such that $U\subset\bigcup_{i}C_i$. For all $n\in\mathbf{N}^*$, we have 
\begin{equation*}
U\subset\bigcup_{1\leq i\leq n}C_i \cup \left[U\setminus\bigcup_{1\leq i\leq n}C_i\right].
\end{equation*}
Then, (\ref{psiunionC}) implies
\begin{align}\label{majpsiC}
\psi(U)\leq\sum_{i=1}^{\infty}\psi(C_i) + 
\psi\left(U\setminus\bigcup_{1\leq i\leq n}C_i\right).
\end{align}
Using $L^2$-monotone outer continuity of $X$ and proposition 1.4.8 in \cite{Ivanoff}, we have
\begin{equation}\label{convpsiC}
\lim_{n\rightarrow\infty}\psi\left(U\setminus\bigcup_{1\leq i\leq n}C_i\right)
=0
\end{equation}
Thus, (\ref{majpsiC}) and (\ref{convpsiC}) imply that for all sequence $\left(C_i\right)_{i\in\mathbf{N}}$ in $\mathcal{C}$ such that $U\subset\bigcup_{i}C_i$,
\begin{equation*}
\psi(U)\leq\sum_{i=1}^{\infty}\psi(C_i)
\end{equation*}
and then, by definition of $m(U)$
\begin{equation}\label{ext2}
\psi(U)\leq m(U).
\end{equation}

\end{itemize}
Equality (\ref{ext}) follows from (\ref{ext1}) and (\ref{ext2}).

From (\ref{defpsi}) and (\ref{ext}), the Borel measure $m$ defined by (\ref{defm}) satisfies
\begin{equation*}
\forall U\in\mathcal{A};\quad
E\left[X_U\right]^2=\psi(U)^{2H}=m(U)^{2H}.
\end{equation*}
\vspace{10pt}

Consider a set-indexed fractional Brownian motion $Y$, defined by
\begin{equation*}
\forall U,V\in\mathcal{A};\quad
E\left[Y_U Y_V\right]=\frac{1}{2} \left[ m(U)^{2H} + m(V)^{2H} - m(U\bigtriangleup V)^{2H} \right]
\end{equation*}
According to proposition 6.4 in \cite{ehem}, projections of $Y$ on any elementary flow $f:[a,b]\rightarrow\mathcal{A}$ is a time-change one-parameter fractional Brownian motion, i. e. such that
\begin{align*}
\forall s,t\in [a,b];\quad
E\left[Y^f_t - Y^f_s\right]^2 &=|m\left[f(t)\right] - m\left[f(s)\right]|^{2H}\\
&=|\theta_f(t) - \theta_f(s)|^{2H},
\end{align*}
where the projection $Y^f$ is defined by $Y^f_t=Y_{f(t)}$, for all $t$.\\
Then, the projections of the set-indexed processes $X$ and $Y$ on any elementary flow have the same distribution. By additivity, this fact holds also on any simple flow. Thus, lemma \ref{lemprojflow} implies $X$ and $Y$ have the same law.
\fin

As a corollary, we get an integral representation.

\begin{corollary}[{\bf Integral Representation}] \label{intcorr}
Let $X=\left\{X_U;\;U\in\mathcal{A}\right\}$ be a $L^2$ outer-continuous set-indexed process. 
Then, $X$ is a sifBm if and only if
for any $U\in\mathcal A$, there exist $f\in F_U^e$ and a Brownian motion $W_f$ such that
\begin{equation}\label{intrep}
X_U=\int_{\mathbf R} (|m(U)-u|^{H-1/2}-|u|^{H-1/2})W_f(du)
\end{equation}
where $H\in[0,1/2)$.
\end{corollary}

\proof
The implication is obvious. Let us prove the converse.\\  
Let $U\in\mathcal A$,
$\forall f\in F_U^e,\ \exists\theta_f: \theta_f(t)=m(f(t))=m(U)$. 
Then
$$X_U=\sifbm^H_{f(t)}=(\sifbm^H)_t^f=\int_{\mathbf{R}}
(|\theta_f(t)-u|^{H-1/2}-|u|^{H-1/2})W_f(du),$$
and the result follows.
\fin

\begin{remark} 
\begin{itemize}
\item If $H=1/2,$ formula (\ref{intrep}) does not hold, but if we decompose $\mathbf R$
into negative and positive parts, the formula can be also interpreted for $H=1/2$.
\item As $W_f$ depends on the flow $f$, expression (\ref{intrep}) can not provide an integral representation of the whole set-indexed process $\sifbm^H$, but only of its projection on a flow.
\end{itemize}
\end{remark}

\bibliographystyle{plain}
\bibliography{style}

\end{document}